\documentclass[12pt]{amsart}

 \newtheorem{proposition}{Proposition}[section]
 \newtheorem{lemma}[proposition]{Lemma}
 \newtheorem{corollary}[proposition]{Corollary}
\newtheorem{theorem}[proposition]{Theorem}

 \theoremstyle{definition}
 \newtheorem{definition}[proposition]{Definition}
 \newtheorem{example}[proposition]{Example}

 \theoremstyle{remark}
 \newtheorem{remark}[proposition]{Remark}

 \newcommand{\thlabel}[1]{\label{th:#1}}
 \newcommand{\thref}[1]{Theorem~\ref{th:#1}}
 \newcommand{\selabel}[1]{\label{se:#1}}
 \newcommand{\seref}[1]{Section~\ref{se:#1}}
 \newcommand{\lelabel}[1]{\label{le:#1}}
 \newcommand{\leref}[1]{Lemma~\ref{le:#1}}
 \newcommand{\prlabel}[1]{\label{pr:#1}}
 \newcommand{\prref}[1]{Proposition~\ref{pr:#1}}
 \newcommand{\colabel}[1]{\label{co:#1}}
 \newcommand{\coref}[1]{Corollary~\ref{co:#1}}
 \newcommand{\relabel}[1]{\label{re:#1}}
 
 \newcommand{\exlabel}[1]{\label{ex:#1}}
 
 \newcommand{\delabel}[1]{\label{de:#1}}
 
 \newcommand{\eqlabel}[1]{\label{eq:#1}}
 \newcommand{\equref}[1]{(\ref{eq:#1})}

 \newcommand{\Hom}{\rm{Hom}}
 \newcommand{\Ext}{\rm{Ext}}

 \newcommand{\Ker}{\rm{Ker}\,}

 \newcommand{\Mod}{\rm{mod}}
 \newcommand{\End}{\rm{End}\,}

 \newcommand{\codim}{\rm{codim}}

 \def\text#1{{\rm {\rm #1}}}

 \def\F{{\mathcal F}}
 \def\T{{\mathcal T}}

 \def\l{\lambda}
 \def\L{\Lambda}
 \def\d{\Delta}
 \def\D{\Delta}
 \def\b{\nabla}

\begin{document}
\title[Good filtration dimensions]{On good filtration dimensions
for standardly stratified algebras} \dedicatory{To Professor
Vlastimil Dlab on the occasion of his 70th birthday.}
\author{Bin Zhu}
\address{Department of Mathematical Sciences, Tsinghua University,
100084 Beijing, China}
 \email{bzhu@math.tsinghua.edu.cn}
 \author{S. Caenepeel}
 \address{Faculty of Applied Sciences,
 Vrije Universiteit Brussel, VUB, B-1050 Brussels, Belgium}
 \email{scaenepe@vub.ac.be}
 \thanks{Research supported by the bilateral project BIL99/43 ``New
 computational, geometric and algebraic methods applied to quantum groups and
 differential operators" of the Flemish and Chinese governments.}
 \urladdr{http://homepages.vub.ac.be/~scaenepe/}
 \begin{abstract}
 $\b$-good filtration dimensions of modules and of algebras are introduced
 by Parker
  for quasi-hereditary algebras. These concepts are now generalized to the
 setting of standardly stratified
 algebras. Let $A$ be a standardly stratified algebra. The $\overline{\b
 }$-good filtration dimension of $A$
 is proved to be the projective dimension of the characteristic module of
 $A$. Several characterizations of
 $\overline{\b }$-good filtration dimensions and $\overline{\d }$-good
 filtration dimensions are given for
 properly stratified algebras. Finally we give an application of these
 results to
 the global dimensions of quasi-hereditary algebras with exact Borel
subalgebra.
 \end{abstract}
 \maketitle

 \section*{Introduction}\selabel{0}
 As generalizations of quasi-hereditary algebras, properly stratified algebras
 and standardly stratifed algebras
  have been introduced by Cline, Parshall, Scott \cite{3} and Dlab \cite{5}.
 They appear in the work of Futorny, K\"onig and
 Mazorchuk on a generalization of the category $O$ \cite{6}. Analogously to
 quasi-hereditary algebras \cite{2}, \cite{4}, \cite{14},  given a standardly
 stratified  algebra $(A,\L)$,  of central importance
  are the modules filted by respectively standard modules, costandard modules,
  proper standard modules or proper costandard
  modules (the precise meaning will be given in \seref{1})
 \cite{1}, \cite{12}.\\
 Recently, in order to calculate the global dimension of the Schur algebra for
 $GL_2$ and $GL_3$,
  Parker \cite{11} introduced the notion of $\b$- (or $\D$-)good
 filtration dimension for
  a quasi-hereditary algebra. The aim of this note is to calculate
 these dimensions for
 standardly stratified algebras and properly stratified algebras.\\
 In \seref{1}, we recall some definitions and results which will be needed
later
 on.
  Various characterizations of $\overline{\b}$-good filtration dimensions of
 standardly stratified
 algebras and of properly stratified algebra are given in \seref{2}; the
 finitistic dimension for
  a properly stratified algebra
  is proved to have an upper bound by the sum of projective dimension and
 injective dimension of its
  characteristic modules.  In \seref{3}, the results of \seref{2} will be
 applied
 to
  quasi-hereditary algebras with a duality and with an exact Borel subalgebra
 in the sense of K\"onig \cite{10}. Their global
 dimensions have an upper bound by two times of that of their exact Borel
 subalgebras.
 This class of algebras includes some algebras arising from Lie theory.

 \section{Preliminary results}\selabel{1}
 Let $R$ be a commutative Artin ring and $A$ a basic Artin algebra
 over $R$.  We will consider finitely generated left $A$-modules;
 maps between $A$-modules will be written on the righthand side of
 the argument, and the composition of maps $f:\ M_1\rightarrow
 M_2,\ g:\ M_2\rightarrow M_3$ will be denoted by $fg$. The
 category of left $A$-modules will be denoted by $A$-mod. All
 subcategories considered will be full and closed under
 isomorphisms.

 Given a class $\Theta$ of $A$-modules, we
 denote by $\F(\Theta)$ the full subcategory of all $A$-modules which
have a $\Theta$-filtration, that is, a filtration
 $$0=M_t\subset M_{t-1}\subset \cdots \subset M_1\subset M_0=M$$
 such that each factor $M_{i-1}/M_i$ is isomorphic to an object of $\Theta$ for
 $1\leq i\leq t$. The modules
 in $\F(\Theta)$ are called $\Theta$-good modules, and the category
 $\F(\Theta)$ is
 called the $\Theta$-good module category.

 In the following, $(A,\le)$ will denote the algebra $A$ together with a fixed
 ordering
 on a complete set $\{ e_1, \cdots , e_n\}$ of primitive orthogonal idempotents
 (given by
 the natural ordering of indices). For $1\le i\le n$ let
 $E(i) $  be the simple $A$-module, which is the simple top of the
 indecomposable
 projective $P(i)=Ae_i$. The standard module $\D (i)$ is by definition the
 maximal factor
 module of $P(i)$ without composition factors $E(j)$ with $j>i$.
 $\overline{\D (i)}$ will be the notation for
 proper standard
 module, which is the maximal factor module of $\D(i)$
  such that the multiplicity condition
 $$[\overline{\D (i)}:E(i)]=1$$
 holds.

 Dually for $1\le \l\le n$, we have costandard modules $\b (\l)$ and
 proper costandard modules $\overline{\b (\l)}$.

 Let $\Delta$ be the
 full subcategory consisting of all $\Delta(\l), \ \l\in \L,$
  and $\Delta_{<\l}$ the full subcategory consisting of all
 $\Delta(\delta), \ \delta <\l$.
 In a similar way, we introduce $\b$ and
 $\b_{<\l} $, and so on.

 The pair $(A,\le )$ is called standardly stratified if ${}_AA\in
 \F(\D)$ (compare \cite{1}, \cite{12}). $(A,\le )$ is called properly
standardly
 stratified if ${}_AA\in \F(\D)$ and ${}_AA\in \F(\overline{\D})$
 (cf. \cite{5}). Note
 that these properties generalize the concept of quasi-hereditary
 algebras where we require the additional condition that the standard
 modules are
 Schur modules.\\

 Let $(A,\le)$  be a standardly stratified algebra.
 It was proved  in \cite{1}, \cite{12} that  $\F(\D)$ and  $\F(\b)$ are
 functorially finite in $A$-mod, which means that they are
 at the same time covariantly
 and contravariantly finite in $A$-mod.  A full subcategory $\T$ of
 $A$-mod is called contravariantly finite in $A$-mod if for any
 $A$-module $M$ there is a module $M_1\in\T$ and a morphism
 $f:\ M_1\rightarrow
 M$ such that the restriction of  $\Hom_A(-,f)$  to  $\T$  is surjective.
 Such a
 morphism $f$  is called a right  $\T$-approximation of $M$.
  A  right $\T$-approximation $f:\ M_1\rightarrow M$ of $M$ is called
minimal if
 the restriction of $f$ to any non-zero direct summand of $M_1$ is non-zero.
 The
 covariant finiteness of $\T$, left $ \T$-approximations of $M$
 and the minimal left $\T$-approximation of $M$ can be defined
 using duality arguments (compare \cite{1}, \cite{12}).\\

 It was also proved in \cite{1}, \cite{12} that there is a unique basic module
 ${}_{A }T$ such that add$({}_{A}T)=\F(\D)\sqcap\F(\overline{\b })$.
 Such a module ${}_{A}T $ is a generalized tilting $A$-module, and is
 called the characteristic module of $A$.  The endomorphism ring
 of ${}_{A}T$ is again a standardly stratified algebra. These
 properties of standardly stratified algebras are summarized in the
 next Lemma.

 \begin{lemma}\lelabel{1.1}
 Let $(A,\le )$  be a standardly stratified
 algebra. Then the following statements hold.
 \begin{enumerate}
 \item $\F(\D)$ is a functorially finite and resolving subcategory
 \item $\F(\overline{\b})$ is a covariantly finite and coresolving subcategory
   of $A$-\Mod.
 \item $\F(\D)=\{ X\in
 A\hbox{-}\Mod~|~\Ext^1(X,\F(\overline{\b}))=0\}$.
 \item $\F(\overline{\b})=\{Y\in
 A\hbox{-}\Mod~|~\Ext^1(\F(\D),Y)=0\}$.
 \item There exists a tilting module ${}_AT$ with
 add$({}_{A}T)=\F(\D)\sqcap\F(\overline{\b })$.
 \end{enumerate}
 \end{lemma}

 It follows from \leref{1.1} that there exists a finite
 $\F(\overline{\b})$-coresolution
 \begin{equation}\eqlabel{1.1.1}
 0\rightarrow X\rightarrow M_0\rightarrow \cdots \rightarrow M_d\rightarrow 0
 \end{equation}
 with $M_i\in\F(\overline{\b})$, for all $X\in A\hbox{-}\Mod$.

 \begin{definition}\delabel{1.1} \cite{11}
 Let $A$ be a standardly stratified algebra. The
 $\overline{\b}$-good filtration dimension of $X$ is the smallest
 number $d$ for which we have an $\F(\overline{\b})$-coresolution
\equref{1.1.1}
 with $M_i\in\F(\overline{\b})$. We then write
 $$\overline{\b}\hbox{-}{\rm gfd}(X)=d$$
 \end{definition}

 Recall the following result from \cite{Friedlander}, \cite{11}:

 \begin{lemma}\lelabel{1.2}
 $\overline{\b}\hbox{-}{\rm gfd}(X)=d$ if and only if
 $\Ext^i_A(\D(\l), X)=0$ for all $i>d$ and all $\l\in \L,$ but there
 exists  $\l\in \L$ such that $\Ext^d_A(\D(\l), X)\not=0$.
 \end{lemma}

 \begin{remark}\relabel{1.1}
 Using the duality principle, we can introduce the dual notions of
 $\overline{\D}$-good module filtration and $\overline{\D}$-good
 filtration dimension of an $A$-module $X$, denoted by
 $\overline{\D}\hbox{-}{\rm gfd}(X)$.
 \end{remark}

 From \cite{11}, we recall the following definition:

 \begin{definition}\delabel{1.2}
 Let $(A,\le)$ be a standardly stratified
 algebra.
 $$d={\rm sup}\{ \overline{\b}\hbox{-}{\rm gfd}(X)~|~
  X\in A\hbox{-}\Mod\}=:\overline{\b}\hbox{-}{\rm gfd}(A)$$
 is called the $\overline{\b}$-good filtration dimension of $A$.
 In a similar way
 $$d={\rm sup}\{ \overline{\D}\hbox{-}{\rm gfd}(X)~|~
  X\in A\hbox{-}\Mod\}=:\overline{\D}\hbox{-}{\rm gfd}(A)$$
 is called the $\overline{\D}$-good filtration dimension of $A$.
 \end{definition}

 The $\overline{\b}$-good filtration dimension
 of $A$ considered as a left $A$-module will be denoted
 denoted by $\overline{\b}\hbox{-}{\rm gfd}({}_AA)$. A similar notation
 will be used for the $\overline{\D}$-good filtration dimension.

 \section{Main result}\selabel{2}
 Throughout this Section, $A$ will be a standardly stratified algebra
  with poset $(\L,\le)$ and $T=\oplus_{\l\in \L}T(\l)$ is the
  characteristic module of $A$.

 \begin{proposition}\prlabel{2.1}
 Let $X\in A\hbox{-}\Mod$. Then
 $\overline{\b}\hbox{-}{\rm gfd}(X)=d$ if and only if $\Ext^i_A(T,X)=0$, for
 any $i>d$, but there exists $\l\in \L,$ such that
 $\Ext^d_A(T(\l),X)\neq 0$.
 \end{proposition}

 \begin{proof}
 Suppose $\overline{\b}\hbox{-}{\rm gfd}(X)=d$. It follows
 from \leref{1.2} that  $\Ext^i_A(T,X)=0$, for any $i>d$. There
 exists $\l\in \L$ such that $\Ext^d_A(\D(\l),X)\neq 0$, and let
 $\l_{0}\in \L$ be minimal with respect to this property.
 Applying $\Hom_A(- ,X)$ to the exact sequence
 $$0\rightarrow
 \D(\l_0)\rightarrow T(\l_0)\rightarrow M(\l_0)\rightarrow 0$$
  where $M(\l_0)\in\F(\D_{<\l_0})$,
 we obtain an exact sequence:
 \begin{equation}\eqlabel{2.1.1}
 \Ext^d(M(\l_0), X)\rightarrow \Ext^d(T(\l_0), X)
 \rightarrow \Ext^d(\D(\l_0),X)\rightarrow
 \Ext^{d+1}(M(\l_0), X)
 \end{equation}
 Therefore $\Ext^d(T(\l_0),X) \simeq \Ext^d(\D(\l_0),X)\neq 0$.\\
 To prove the converse, we use induction on $\l$ to prove that
 $\Ext^i(\D(\l),X)=0$ for all $\l$ and $i>d$.
  If $\l$ is minimal in $\L$, then $\D(\l)=T(\l)$, and we are done.
 Take a non-minimal $\l\in\L$; applying $\Hom_A(-,X)$ to the exact sequence
 $$0\rightarrow \D(\l)\rightarrow
 T(\l)\rightarrow M(\l)\rightarrow 0$$
 where $M(\l)\in\F(\D_{<\l})$, we obtain an exact sequence
 \begin{equation}\eqlabel{2.1.2}
 \Ext^i(M(\l), X)\rightarrow \Ext^i(T(\l), X)
 \rightarrow \Ext^i(\D(\l),X)\rightarrow
 \Ext^{i+1}(M(\l), X)
 \end{equation}
  It follows that $\Ext^i(\D(\l),X)=0$ for all $i>d$, and
 we have that $\overline{\b}\hbox{-}{\rm gfd}(X)\le d$.\\
 Let $\l$ be the minimal weight such that
    $\Ext^d(T(\l),X)\not=0$. It follows from the exact sequence \equref{2.1.2}
 (with i=d) that
 $$\Ext^d(\D(\l), X)
 \simeq \Ext^d(T(\l),X)\not=0$$
 Then $\overline{\b}\hbox{-}{\rm gfd}(X)=
 d$, finishing our proof.
 \end{proof}

 As a consequence, we find

 \begin{proposition}\prlabel{2.2}
 $\overline{\b}\hbox{-}{\rm gfd}(A)=s$ if and only
 if ${\rm proj.dim}_A T=s$.
 \end{proposition}

 \begin{proof}
 $\overline{\b}\hbox{-}{\rm gfd}(A)=s$ if and only if
 $\overline{\b}\hbox{-}{\rm gfd}(X)\le s$ for all
 $X\in A\hbox{-}\Mod$, and there is an $A$-module $M$ such that
   $\overline{\b}\hbox{-}{\rm gfd}(M)=s$.
 By \prref{2.1}, this is equivalent to
 $\Ext^j(T,X)=0$, for all $j>s$ and the existence of an $A$-module $N$ such
 that $\Ext^s(T,X)=0$. The last condition is equivalent to
 ${\rm proj.dim}_A T=s$, finishing the proof.
 \end{proof}

 \begin{proposition}\prlabel{2.3}
 ${\rm proj.dim}_A T=d$ if and only if
 $\Ext^i(T,A)=0$ for all $i>d$, and $\Ext^d(T,A)\neq 0.$
 \end{proposition}

 \begin{proof}
 Assume that ${\rm proj.dim}_A T=d$. It is easy to check that
 $\Ext^d(T,A)\neq 0$.\\
 The proof of the converse consists of two steps.
 We first prove that for any $\l$ and any $j>d$,
 $\Ext^j(T, \D(\l))=0$. If $\l$ is maximal in $\L$, then
 $\D(\l)=T(\l)$, and then $\Ext^j(T, \D(\l))=0$ for any $j>d$.
 If
  $\l$ is not maximal, then we have an exact sequence
 $$0\rightarrow
 U(\l)\rightarrow P(\l)\rightarrow  \D(\l)\rightarrow 0$$
where $U(\l)\in
 \F(\D_{>\l})$.
 Applying $\Hom(T,-)$, we find
 the exact sequence:
 $$\Ext^j(T,P(\l))\rightarrow
 \Ext^j(T, \D(\l)) \rightarrow \Ext^{j+1}(T, U(\l))$$
 and it follows that
 $$\Ext^j(T, \D(\l))=0$$
 for any $j>d$.\\
  Secondly, for any $M\in A\hbox{-}\Mod$ we have that $\Ext^j(T,
  M)=0$, for any $j>d$.  Let $ f:\ X\rightarrow M$ be the minimal
  $\F(\D)-$approximation of $M$. Then $f$ is surjective and $K=\Ker f\in
  \F(\b),$ and we have the exact sequence
 $$0\rightarrow K\rightarrow
 X\rightarrow M\rightarrow 0$$
 Applying $\Hom(T,-)$, we find an exact sequence
 $$\Ext^j(T,X))\rightarrow \Ext^j(T, M) \rightarrow
 \Ext^{j+1}(T, K)$$
  From the fact that $\Ext^{j+1}(T,
 K)=0$, it follows that $\Ext^j(T, M)=0$ for all $j>d$, and therefore
 ${\rm proj.dim}_AT=d$.
 \end{proof}

 Recall from \cite{8}, we recall the notion of $T$-codimension of
 an $A$-module $X$: it is the smallest number $s$ such that we have
 an exact sequence
 $$ 0\rightarrow X\rightarrow T_0\rightarrow \cdots \rightarrow
 T_s\rightarrow 0$$
 with $T_i\in {\rm add}T$ for all $0\le i\le s$. We then write
 $$T\hbox{-}\codim(X)=s$$
 In a similar way, we introduce $T\hbox{-}\dim(X)$, see \cite{15}.

 \begin{proposition}\prlabel{2.4}
 ${\rm proj.dim}(T)= T\hbox{-}\codim({}_AA)$.
 \end{proposition}

 \begin{proof}
 It follows from \cite[Lemma 2.2]{8} that
 $T\hbox{-}\codim({}_AA)\leq {\rm proj.dim}(T)$.
 Let $T\hbox{-}\codim({}_AA)=d$. Then we have a $T$-coresolution of
 $A$ of length $d$, namely
 $$ 0\rightarrow A\rightarrow T_0\rightarrow \cdots \rightarrow
 T_d\rightarrow 0$$
 with $T_i\in {\rm add}T$ for $0\le
 i\le d$.  It follows from \cite{8} that any indecomposable direct
 summand of $T$ appears in the coresolution, and therefore we have that
 $\Ext^j(T,A)=0$ for any $j>d$. It follows from \prref{2.3}.
 that ${\rm proj.dim}(T)\le d$, as needed.
 \end{proof}

 Combining everything together, we find the following descriptions of the
 $\overline{\b}$-good filtration dimension of a standard stratified algebra
$A$.

 \begin{theorem}\thlabel{2.5}
 Let $(A,\le)$ be a standard stratified
 algebra, $T$ the characteristic (full tilting) module of $A$ and
 $d$ a non-negative integer. Then the following are equivalent:
 \begin{enumerate}
 \item $\overline{\b}\hbox{-}{\rm gfd}(A)=d$,
 \item ${\rm proj.dim}(T)=d,$
 \item $\overline{\b}\hbox{-}{\rm gfd}(_AA)=d$,
 \item $T\hbox{-}\rm codim(_AA)=d$.
 \end{enumerate}
 \end{theorem}

 \begin{proof}
 The equivalence of 1) and 2) is just \prref{2.2},
 the equivalence of 2) and 4) is \prref{2.3}, and the
  equivalence of 2) and 3) is \prref{2.4}.
 \end{proof}

 \begin{remark}\relabel{2.1}
 Let $A$ be a quasi-hereditary algebra. According to the
 Definitions in \seref{1}, the good filtration dimension of $A$
 as an algebra, and $A$ considered as an $A$-module, could be
 different (see also the remark following Definition 2.3
 in \cite{11}).
 The equivalence of 1) and 4) in \thref{2.5} tells us that they are equal:
 \begin{equation}
 \eqlabel{2.1.3}
 \overline{\b}\hbox{-}{\rm gfd}(A)=
  \overline{\b}\hbox{-}{\rm gfd}(_AA)
 \end{equation}
 \end{remark}

 In the following we deal with properly stratified
  algebras; for basic properties of these algebras we refer to
 \cite{5}. We are indebted to the referee for pointing out to us
 that \thref{2.6} also follows from \cite{1}, together
 with Dlab's result
 (see \cite{5}) that the oppposite algebra of a properly stratified algebra
 is standardly stratified. We kept a short proof for the sake of
 completeness.

 \begin{theorem}\thlabel{2.6}
 Let $(A,\le)$ be a properly stratified
  algebra. Then there exist a tilting module $T$ and a cotilting module $S$
 such that
   ${\rm add}T=\F(\D)\sqcap\F(\overline{\b })$ and
 ${\rm add}S=\F(\b)\sqcap\F(\overline{\D })$.
 \end{theorem}

 \begin{proof}
 Let $(A,\le)$ be a properly stratified
  algebra. Then $(A,\le)$ is standardly
 stratified, and then there exists a tilting module $T$ with
   ${\rm add}T=\F(\D)\sqcap\F(\overline{\b })$. Since ${}_AA\in
 \F(\overline{\D }),$  we have $D(A_A)\in
  \F(\b)$. The dual version of \leref{1.1} or \cite[Theorem 2.1]{1}
 tells us that
 $$\F(\overline{\D })={}^{\bot}\F(\b){\rm  ~~ and  ~~}
 \F(\b)={}^{\bot}\F(\overline{\D})$$
 By the dual version of \cite[Theorem 2.1]{1}, we have
 a basic cotilting $A$-module $S$ with
 ${\rm add}S=F(\b)\sqcap\F(\overline{\D})$.
 \end{proof}

 \begin{corollary}\colabel{2.7}
 Let $(A,\le)$ be a properly stratified algebra and
 take $T$ and  $S$ as in \thref{2.6}. Then
 \begin{eqnarray}
 \overline{\b}\hbox{-}{\rm gfd}(A)&=&{\rm proj.dim}(T)=
 T\hbox{-}{\rm codim}({}_AA)=\overline{\b}\hbox{-}{\rm gfd}({}_AA)
 \eqlabel{2.7.1}\\
 \overline{\D}\hbox{-}{\rm gfd}(A)&=&{\rm inj.dim}(S)=
 S\hbox{-}{\rm dim}(D(A_A))=\overline{\D}\hbox{-}{\rm gfd}(D(A_A))
 \eqlabel{2.7.2}
 \end{eqnarray}
 \end{corollary}

 \begin{proof}
 \equref{2.7.1} follows from \thref{2.5}, and \equref{2.7.2} is the dual
 version of it.
 \end{proof}

 Applying \coref{2.7} to quasi-hereditary algebras, we recover the
 various descriptions of good filtration dimensions introduced in
 \cite{11}.

 \begin{corollary}\colabel{2.8}
 Let $(A,\le)$ be a quasi-hereditary
 algebra and $T$ the characteristic module of $A$. Then
 \begin{eqnarray}
 \b\hbox{-}{\rm gfd}(A)&=&{\rm proj.dim}(T)
 =T\hbox{-}{\rm codim}({}_AA)=\b\hbox{-}{\rm gfd}({}_AA)\eqlabel{2.8.1}\\
 \D\hbox{-}{\rm gfd}(A)&=&{\rm inj.dim}(T)=T\hbox{-}{\rm
 dim}(D(A_A)) =\D\hbox{-}{\rm gfd}(D(A_A))\eqlabel{2.8.2b}
 \end{eqnarray}
 \end{corollary}

 \begin{proof}
 If $(A,\le)$ be quasi-hereditary, then
 $\overline{\D}(i)=\D(i)$ and $\overline{\b}(i)=\b(i)$, and therefore
 $S=T$. All the assertions follow.
 \end{proof}

 \begin{theorem}\thlabel{2.9}
 Let $(A,\le)$ be a properly stratified algebra with $T=S$,
 where $T$ and $S$ are as in \thref{2.6}. Then the
 sum of the projective and injective dimension of $T$ is finite,
 and is an upper bound for the finitistic dimension of $A$.
 \end{theorem}

 \begin{proof}
 Let $(A,\le)$ be a properly stratified
 algebra with $T=S$. Then ${\rm proj.dim}(T)=t <\infty$ and
 ${\rm inj.dim}(T)=s<\infty$. Let $X$ be any left $A$-module with finite
 projective dimension and $\Omega ^s(X)$ the $s$-th syzygy module of $X$. It
 follows that $\Ext^{s+i}(X,T)=0$ for any
 $i>0$. Therefore $\Ext^{i}(\Omega^s(X),T)=0$ for any $i>0$. One has
 that $\Omega^s(X)\in {}^{\bot}T$. It follows from
 ${\rm proj.dim}(\Omega^s(X))<\infty$ that $\Omega^s(X)\in \F(\D)$ (see
 \cite{11}).
  Therefore we have that ${\rm proj.dim}(\Omega^s(X))\le t$, and
 ${\rm proj.dim}(X)\le t+s$. It follows that the
 finitistic dimension of $A$ is at most $s+t$.
 \end{proof}

 Note that quasi-hereditary algebras are examples of properly
  stratified algebras satisfying the condition in \thref{2.9}.
 A non-quasi-hereditary example is the following: let
 $A$ be the quiver algebra of the quiver
 consisting of one vertex and one arrow $\alpha$, with relation
 $\alpha ^2=0$. Then $\F(\D)=\F(\b)={\rm add}(A)$. It is easy to
 see that $A$ is a properly stratified algebra with $T=S$.\\

 Applying \thref{2.9} to a quasi-hereditary algebra $A$, we find that
 the global dimension of a quasi-hereditary algebra is bounded by the sum of
 the projective and injective dimension
  of their characteristic module $T$. We remark at this point that the
same formula holds for the global dimension of the endomorphism
algebra
of a tilting module $T$ over an Artin algebra, see \cite{7}.\\

 As a direct consequence of \thref{2.9}, we have

 \begin{corollary}\colabel{2.10}
 Let $A$ be a quasi-hereditary algebra and
 $T$ the characteristic module of $A$. Then
 $$\max\{{\rm proj.dim}(T),{\rm inj.dim}(T)\}\le
 {\rm gl.dim}(A)\le {\rm proj.dim}(T)+{\rm inj.dim}(T)$$
 \end{corollary}

 \begin{example}\exlabel{2.10}
 There exist quasi-hereditary algebras such that the right
 inequality in \coref{2.10} is strict.
 Let $A$ be the path algebra given by
 \begin{center}
 \setlength{\unitlength}{0.71cm}
 \begin{picture}(8,1)
 \multiput(0.5,0)(2,0){2}{\vector(1,0){1.7}}
 \multiput(0.2,0)(2.1,0){3}{\circle{0.1}} \put(0.2,0.2){$2$}
 \put(2.4,0.2){$1$} \put(4.5,0.2){$3$}
 \end{picture}
 \end{center}
 For the ordering $\ \L =  \{ \  1<2<3\ \},\ \  (A,\L) \ $ is a
 quasi-hereditary algebra with characteristic module $T=E(1)\oplus
 P(1)\oplus P(2)$. Then ${\rm gl.dim}(A)<
 {\rm proj.dim}(T)+{\rm inj.dim}(T)$.
 \end{example}

 Let $A'={\rm End}_AT$. $A'$ is also a quasi-hereditary algebra
 with respect to the opposite ordering of $\L$, and is called the
 Ringel dual of $A$. $T'={\rm Hom}_A(T,D(A))$ is the characteristic
 module of $A'$, and the endomorphism ring of ${}_{A'}T' $ is Morita
 equivalent to $A$ as a quasi-hereditary algebra. Applying
 \thref{2.9} to the Ringel dual of $A$, we obtain a description of
 its global dimension.

 \begin{corollary}\colabel{2.11}
 Let $A$ be a quasi-hereditary algebra, $T$
 its characteristic module, and $A'$ its Ringel dual.
 Then
 \begin{eqnarray*}
 \max\{{\rm proj.dim}({}_{A}T),{\rm inj.dim}({}_{A}T)\}&\le&
 {\rm gl.dim}(A')\\
 &\le& {\rm proj.dim}({}_{A}T)+{\rm inj.dim}({}_{A}T)
 \end{eqnarray*}
 \end{corollary}

 \begin{proof}
 We assume that $T$ is the basic characteristic module, that is, it
 contains exactly $|\L|$ indecomposables.  Let ${\rm
 proj.dim}({}_AT)=s$, ${\rm inj.dim}({}_AT)$ $=t$. Then from \equref{2.8.2b} in
 \coref{2.8} and the fact that the endomorphism ring
   $\End_{A'}T'$ is isomorphic to $A$ as a quasi-hereditary algebra \cite{14},
   we have that ${\rm proj.dim}({}_{A'}T')=t$ and
 ${\rm proj.dim}(_{A'}T')=s$. By \coref{2.10}, we have that
 $${\rm gl.dim}(A')\le{\rm proj.dim}({}_{A}T)+{\rm inj.dim}({}_{A}T)$$
 \end{proof}

 \begin{remark}\relabel{2.2}
 It would be interesting to know for which
 quasi-hereditary algebra $A$ we have the equality
 \begin{equation}\eqlabel{2.2.1}
 {\rm gl.dim}(A)=
 {\rm proj.dim}({}_AT)+{\rm inj.dim}({}_AT)
 \end{equation}
 It was proved in \cite{Totaro} and \cite[Sec. 4.8]{Donkin}
 \equref{2.2.1} holds for Schur algebras $S(n,r)$, with $p$
 arbitrary and $n\geq r$. Parker \cite{Parker2} recently proved that
 \equref{2.2.1} holds for Schur algebras with $n=2,3$ and $p$
 arbitrary, or $p>n$ and $r$ arbitrary. In \seref{3}, we will show
 that this equality holds for quasi-hereditary algebras with simple
 modules as standard modules.
 \end{remark}

\section{Applications and Examples}\selabel{3}
 In this section we will apply the results of \seref{2} to some
 special classes of quasi-hereditary algebras, especially to algebras with
 exact Borel subalgebras.
 For the definition of an exact Borel subalgebra of a quasi-hereditary algebra,
 we refer to \cite{10} or \cite{13}. Let $A$ be a quasi-hereditary
 algebra. We say that $A$ has a duality if
 there exists an involutory, contravariant functor $\phi:\
 A\hbox{-}\Mod\rightarrow A\hbox{-}\Mod$ preserving simple
 modules, i.e.
 $\phi (E(i))=E(i)$ for any $i$.
 If $A$ is  a quasi-hereditary algebra $A$ with a duality $\phi$,
 then
 $$\phi (\d(i))=\b(i)~~{\rm and}~~\phi (T(i))=T(i)$$
 for all $i$.

 \begin{lemma}\lelabel{3.1}
 Let $(A,\L)$ be a
 quasi-hereditary algebra with simple modules as $\D-$modules. Then
 $${\rm gl.dim}(A)={\rm proj.dim}(D(A))={\rm inj.dim}({}_AA)$$
 \end{lemma}

 \begin{proof}
 If $A$ is a quasi-hereditary algebra with simple standard modules,
 then $\F(\D)=A\hbox{-}\Mod$. Therefore the
  injective module is the characteristic module of $A$.
 By \coref{2.8} or \thref{2.9}, we have that
 ${\rm gl.dim}(A)={\rm proj.dim}(D(A))$.  $A$ is also
   a quasi-hereditary algebra with respect to the opposite ordering of $\L$.
   For this quasi-hereditary algebra $(A,\L^{op})$,  $\D$-modules are
 indecomposable
   projective modules, and its characteristic module is ${}_AA$.
   Then ${\rm gl.dim}(A)={\rm inj.dim}(A)$.
 \end{proof}

 \begin{lemma}\lelabel{3.2}
 Let $B$ be an exact Borel subalgebra of a
 quasi-hereditary algebra $A$. Then for  any $B-$module ${}_{B}M$,
 the $A-$module $A\otimes_BM$ has a
  $\D_A-$filtation. Moreover
 $${\rm proj.dim}
 (A\otimes_BM)\le {\rm proj.dim}({}_BM)$$
 \end{lemma}

 \begin{proof}
 First, we prove that $A\otimes_BM\in \F(_A\D)$.
  Since $B$ is a Borel subalgebra of $A$, the functor
 $$F=A\otimes-:\ B\hbox{-}\Mod\to A\hbox{-}\Mod$$
 is exact and $F(\D_B(\l))\simeq \D_A(\l)$ (compare \cite{10},
 \cite{13}).
  Let $G:\ A\hbox{-}\Mod\to B\hbox{-}\Mod$ be the restriction
 of scalars functor. We first prove that\\ $G(\b_A(j))\in \F({}_B\b)$
   for all $j$. This follows from the isomorphisms
 \begin{eqnarray*}
 {\rm Ext}^1_B(\D_B(i),G(\b_A(j)))&\simeq& {\rm
 Ext}^1_A(F(\D_B(i)),\b_A(j)))\\
 &\simeq& {\rm Ext}^1_A(\D_A(i),\b_A(j))=0
 \end{eqnarray*}
 for any $i, j\in \L$. Then from the isomorphisms
 $${\rm Ext}^1_A(F(_BM),\b_A(j))\simeq {\rm Ext}^1_B(_BM,G(\b_A(j)))=0$$
 for any
 $ j\in \L$, we have that $F(_BM)\in \F(_A\D)$.
  Since $A\otimes-$ is an
 exact functor, applying
   $A\otimes-$ to a minimal
  projective resolution of ${}_BM$ gives a projective resolution of the
 characteristic
  module $A\otimes_BM$. This means that
 ${\rm proj.dim}(A\otimes _BM)\le {\rm proj.dim}({}_BM)$. This
 finishes the proof.
 \end{proof}

 \begin{theorem}\thlabel{3.3}
 Let $A$ be a quasi-hereditary algebra with a duality and with an exact Borel
 subalgebra $B$. Then
 ${\rm gl.dim}(A)\le 2{\rm gl.dim}(B)$.
 \end{theorem}

 \begin{proof}
 For a quasi-hereditary algebra $A$ with a duality $\phi$, the characteristic
 module $T$ is self-dual, i.e. $\phi (T)\simeq T$. Then
 ${\rm proj.dim}(T)={\rm inj.dim} (T)$, and it follows from Lemmas
 \ref{le:3.1} and \ref{le:3.2} and \coref{2.10}
 that
 ${\rm gl.dim}(A)\le 2{\rm proj.dim}(T)=2\mbox{sup}\{ \mbox{proj.dim}\D_A(\l) | \l\in \L \}
 \le 2\mbox{sup}\{ {\rm proj.dim} \D_B(i) | \l \in \L\}= 2{\rm gl.dim}(B)$.
 \end{proof}

 \begin{example}\exlabel{3.1}
 This example, taken from \cite{10}, shows that it is possible
 that ${\rm gl.dim}(A)= 2{\rm gl.dim}(B)$
 for a Borel subalgebra $B$ of a quasi-hereditary algebra $A$.\\
 Let $A$ be the algebra given by
 $$\begin{array}{c}
 \setlength{\unitlength}{0.7cm}
 \begin{picture}(8,2)
 \put(0.8,1.8){$\alpha$} \put(2.8,1.8){$\gamma$}
 \multiput(0.3,1.5)(2,0){2}{\vector(1,0){1.7}}
 \multiput(4,1)(-2,0){2}{\vector(-1,0){1.7}}
  \put(0.8,0.5){$\beta$}
  \put(2.8,0.5){$\delta$}
 \multiput(0.2,1.2)(2,0){3}{\circle{0.1}} \put(0,0.2){$3$}
 \put(2.1,0.2){$1$} \put(4,0){$2$}
 \end{picture}\end{array}$$
 with relations
 $$\delta\cdot\gamma \cdot \alpha=0,~~~\beta
 \cdot \delta \cdot \gamma=0,~~~\beta \cdot \alpha=0,~~~
 \gamma\cdot\delta=0$$
 Ordering the weights by $1<2<3$, $A$ is a quasi-hereditary
   algebra which has an exact Borel subalgebra $B$ given by
 $$\begin{array}{c}
 \setlength{\unitlength}{0.7cm}
 \begin{picture}(8,2)
 \put(3,1.5){\vector(-1,0){2.5}} \put(3,1.3){\vector(-1,-1){1.3}}
 \put(1.5,0.2){\vector(-1,1){1}}
  \put(1.85,1.7){$\beta$}
  \put(2.4,0.5){$\gamma$}
 \put(0.2,0.4){$\delta \beta$}
 \multiput(0.2,1.4)(3,0){2}{\circle{0.1}}
 \put(1.6,0){\circle{0.1}}
 \put(0.2,1.6){$3$} \put(2.8,1.6){$1$} \put(1.2,0){$2$}
 \end{picture}\end{array}$$
 with relation
 $$\beta\cdot\delta\cdot\gamma =0$$
 It is easy to see that ${\rm gl.dim}(B)=2$ and
 ${\rm gl.dim}(A)=2{\rm gl.dim}(B)=4$.
 \end{example}

 \begin{center}
 {\sc Acknowledgements}
 \end{center}
 The authors would like to thank A. Parker, V. Mazorchuk and
 Changchang Xi for their comments to an earlier version of this
 paper.

 \end{document}